\documentclass[11pt,twoside]{article}
\makeatletter
\@addtoreset{equation}{section}
\makeatother

\def\pref#1{(\ref{#1})}
\def\dsp{\displaystyle}
\def\Frac#1#2{\frac
{
 {\raise.6ex
 \hbox{$\displaystyle#1$}}
}
{
 {\lower.6ex
 \hbox{$\displaystyle#2$}}
 }
}

\input epsf

\def\intp{\int_0^\infty}

\def\bo{{\cal O}}

\def\wt{\widetilde}
\def\tfrac#1#2{{{\lower.6ex
\hbox{$\scriptstyle#1$}}\over 
{\raise.7ex
\hbox{$\scriptstyle#2$}}}}

\begin{document}
 \title{
On the Temporal Order of First-Passage Times in 
One-Dimensional Lattice Random Walks}

\author{J.B. Sanders\\
        FOM Institute for Atomic and Molecular Physics (AMOLF), \\
        Kruislaan 407, 1098 SJ Amsterdam,  The Netherlands\\
        \\
    N.M. Temme\\
        Centrum voor Wiskunde en Informatica (CWI), \\
        Kruislaan 413, 1098 SJ Amsterdam,  The Netherlands\\
    {\ }\\
     { \small e-mail: {\tt
    nicot@cwi.nl}}
    }

 \maketitle 
 
 \begin{abstract} 
 \noindent 
A random walk problem with particles on discrete double infinite linear
grids is discussed.  The model is based on the work of Montroll and others.
A probability connected with the problem is given in the form of integrals
containing modified Bessel functions of the first kind.  By using several
transformations simpler integrals are obtained from which for two and three
particles asymptotic approximations are derived for large values of the
parameters.  Expressions of the probability for $n$ particles are also
derived.

\end{abstract}

\hspace*{\fill}
\begin{minipage}[b]{53ex}
{\tiny{\em
I returned and saw under the sun, that the race is not to the swift,
nor the battle to the strong, neither yet bread to the wise, nor yet
riches to men of understanding, nor yet favour to men of skill; but
time and chance happeneth to them all. }\\
George Orwell, Politics and the English Language, Selected Essays,
Penguin Books,  1957. (The citation is from Ecclesiastes 9:11.)}

\end{minipage}

\vskip 0.8cm \noindent
2000 Mathematics Subject Classification:
%\par\noindent
41A60, 60G50, 33C10.
\par\noindent
Keywords \& Phrases:
%\par\noindent
 random walk,
 asymptotic expansion, 
 modified Bessel function.

\section{Introduction}
\label{sec:int}

The subject of random motion is one on which an enormous amount of
mathematical studies have been made.  We mention in this respect the
classical work of Rayleigh, Smoluchowski, Chandrasekhar, and
countless others \cite{Montroll1979}.  In this paper, we are interested in the
specialization of this general notion to random walk on a periodic
lattice, where a particle makes random jumps between neighbouring
sites of this lattice.  In this respect we refer in particular to
the pioneering work by Montroll and his collaborators which has
provided the inspiration for the present work.

We shall very briefly indicate the method of Montroll's approach,
where throughout this paper we shall limit ourselves to random walks on
one (or more)  linear (1D) lattice chains. We shall also suppose that the jump probabilities
of a random walker to the left and to the right are equal, and hence equal to $p=\frac12$.
Initially, the time is considered to be discrete,
which means that we consider the situation of the particle after a 
discrete number of jumps $n$,
which is equivalent  to allowing the particle to jump once in
every unit of time. Montroll et al. \cite{Montroll1964}, \cite{Montroll1965} now introduce two quantities 
which are of very great importance. These are

\begin{enumerate}
\item
$P_n(\ell)$, the probability that the random walker will be at site $\ell$ after the $n$th jump.
\item
$f_n(\ell)$, the probability that the random walker will be at site $\ell$ after the $n$th jump
{\it for the first time}.
\end{enumerate}
Of course, it is assumed that before the first jump ($n=0$) the 
particle is at the origin ($\ell=0$). 

The function $P_n(\ell)$ satisfies the following equation
\begin{equation}\label{j1}
P_n(\ell)=\tfrac12P_{n-1}(\ell-1)+\tfrac12P_{n-1}(\ell+1).
\end{equation}
(If at epoch $n-1$, the particle is at either $\ell-1$ or $\ell+1$, it will
have a probability $p=\frac12$ to be at $\ell$ at epoch $n$. If it is anywhere else at epoch
$n-1$, its chance of being at $\ell$ one jump later, is zero.)
This equation also shows that  the random walk, as described above, 
is a Markoff process, in that the state $(\ell)$ of the random walker at a given epoch $n$ depends only on that 
at {\it one} moment earlier.

Montroll then introduces a generating function
$U(\ell,z)=\sum_{n=0}^\infty P_n(\ell) z^n$. This function $U(\ell,z)$ is then calculated explicitly, 
from which $P_n$ and various moments over $\ell$ can be calculated. 
For details we refer to \cite{Montroll1964} and \cite{Montroll1965}. We also refer to these papers for the 
treatment of the first passage times $f_n(\ell)$ and the corresponding generating function 
$F(\ell, z)=\sum_{n=0}^\infty f_n(\ell)z^n$. The quantity $f_n(\ell)$ is the probability of 
reaching the site $\ell$ for the first time at the $n$th jump.

For the sake of completeness we give the explicit expressions for $U(\ell,z)$ and $F(\ell,z)$:
\begin{equation}\label{j2}
\begin{array}{l}
\dsp{U(\ell,z)=\left(\frac z2\right)^{\ell}{}_2F_1\left(\frac{\ell+1}{2},
\frac{\ell+2}{2};\ell+1;z^2\right)=
\frac{1}{\sqrt{1-z^2}}\left(\frac{1-\sqrt{1-z^2}}{z}\right)^\ell,}\\
\\
\dsp{F(\ell,z)=\left(\frac z2\right)^{\ell}{}_2F_1\left(\frac{\ell}{2},
\frac{\ell+1}{2};\ell+1;z^2\right)-\delta_{\ell,0}\sqrt{1-z^2}=\frac{U(\ell,z)-\delta_{\ell,0}}{U(0,z)},}
\end{array}
\end{equation}
from which explicit forms of $P_n(\ell)$ and $f_n(\ell)$ follow. 

Montroll et al. \cite{Montroll1965} also present a method of treating the time as a continuous variable. Then we 
introduce as fundamental quantities the following probability densities:
\begin{equation}\label{j3}
\begin{array}{ll}
\dsp{\overline{P}(\ell,t)\,dt:} &\hbox{\rm the probability density for the random walker to {\it be} at $\ell$ }\\
&\hbox{during interval $(t,t+dt)$}. \\ &
\\
\dsp{\overline{F}(\ell,t)\,dt:} &\hbox{\rm the probability density for the random walker to {\it arrive} at $\ell$ }\\
&\hbox{during interval $(t,t+dt)$ for the first time.}
\end{array}
\end{equation}
Jumps are now taken to occur  at random times $t_1,t_2, t_3,\ldots$. This implies the introduction of the 
random variables $T_1=t_1, T_2=t_2-t_1, \ldots T_n=t_n-t_{n-1}$, 
which have the common density $\psi(t)$. For $\psi(t)$ we take the exponential 
density $\psi(t)=\alpha e^{-\alpha t}$, where $\alpha$ is the average number of jumps 
made by the random walker per unit of time.
From this point on we shall concentrate on the first-passage probability
density function, that being the one which we shall need most in 
in the following applications. 

We also introduce the  probability densities
\begin{equation}\label{j4}
\psi_0(t)=\delta(t),\quad \psi_n(t)=\int_0^t \psi(t-\tau)\psi_{n-1}(\tau)\,d\tau, \quad n=1,2,3,\ldots\,.
\end{equation}
The function $\psi_n(t)$ can be interpreted as the probability density that the $n$th jump
of the random walker takes place in the time interval $(t,t+dt)$.
We have
\begin{equation}\label{j4a}
\psi_n(t)=\alpha^n \, e^{-\alpha t}\,\frac{t^{n-1}}{(n-1)!}, \quad n=1,2,3,\ldots\,.
\end{equation}

It can now easily been understood that \cite{Montroll1965}
\begin{equation}\label{j5}
\overline{F}(\ell,t)=\sum_{n=0}^\infty f_n(\ell) \psi_n(t).
\end{equation}
If we use the $\psi_n(t)$  given above and the $f_n(\ell)$ that follow from the second line of \pref{j2},
we obtain 
\begin{equation}\label{j5a}
\overline{F}(\ell,t)=2^{-\ell}e^{-\alpha t} t^{-1}
\sum_{n=0}^\infty \frac{(\ell/2)_n\,(\ell/2+1/2)_n}{n!\,(\ell+1)_n}
\frac{(\alpha t)^{\ell+2n}}{(\ell+2n-1)!}.
\end{equation}
where $(a)_n$ denotes Pochhammer's symbol defined by 
\begin{equation}\label{j5b}
(a)_0=1, \quad (a)_n=a(a+1)\cdots(a+n-1), \quad n=1,2,3,\ldots.
\end{equation}
Comparing the expansion in \pref{j5a} with that of the modified Bessel function of the first kind, 
see \cite[Chapter~9]{Abr}, 
\begin{equation}\label{j5c}
I_\nu(z)=\sum_{n=0}^\infty \frac{(z/2)^{\nu+2n}}{n!\,\Gamma(\nu+n+1)},
\end{equation}
and using the duplication formula of the gamma function
\begin{equation}\label{j5d}
\Gamma(2z)=\frac{2^{2z-1}}{\sqrt\pi}\Gamma(z)\Gamma(z+\tfrac12),
\end{equation}
an explicit form for $\overline{F}(\ell,t)$ is obtained:
\begin{equation}\label{j6}
\overline{F}(\ell,t)=e^{-\alpha t}\frac{\ell }{t} I_{\ell}(\alpha t).
\end{equation}
From \cite[Eq. 11.4.13]{Abr}
it follows that for $\ell\ne0$
\begin{equation}\label{j7}
\intp \overline{F}(\ell,t)\,dt=\ell\intp e^{-\alpha t} I_{\ell}(\alpha t) \frac{dt}{t}=1.
\end{equation}

\begin{figure}
\begin{center}
\epsfxsize=12cm \epsfbox{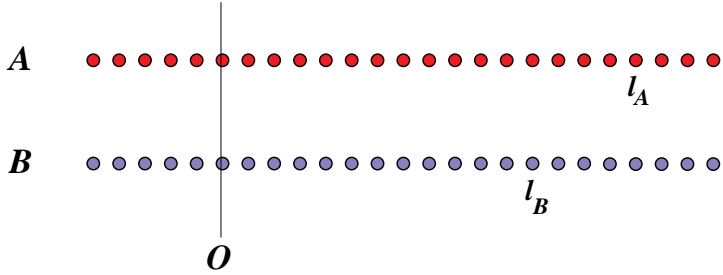}
\end{center}
\caption{\small The random walks for $A$ and $B$. \label{walks1}}
\end{figure}

It is of interest to consider the problem of several simultaneous random walkers on a lattice chain, 
and the behaviour in time of their mutual configuration. It is as an 
introduction to this work that we shall consider two, three, $\ldots$, independent
random walkers on separate lattice chains. 
We begin with two random walkers and consider the 
situation as shown in Figure \ref{walks1}.

\medskip
\noindent
{\bf Remark 1.\quad}
The integers $\ell_A$ and $\ell_B$ may separately assume negative values.
However, to avoid the use of absolute value signs, we 
consider only positive values of $\ell_A$ and $\ell_B$. But all results hold
for negative values when we replace these quantities by their absolute values.

\bigskip

We are interested in finding the probability that particle $A$ arrives at 
$\ell_A$ {\it before} particle $B$ arrives at $\ell_B$.
The solution to this problem is an intermediate 
result for the treatment of a $1D-$grid random walk problem of 
an agglomeration of many particles.

We know that the probability density for $A$ to arrive for the first time at $\ell_A$ 
in the interval $(t,t+dt)$ is 
\begin{equation}\label{j8}
\overline{F}(\ell_A,t)\,dt=\frac{\ell_A}{t}  e^{-\alpha t} 
I_{\ell_A}(\alpha t)\,dt.
\end{equation}
It is now obvious that the desired probability can be written as
\begin{equation}\label{i1}
P(t_{\ell_A}\le t_{\ell_B}) = \int_0^\infty \,dt_A\, \overline{F}(\ell_A,t_A)\  
        \int_{t_A}^\infty dt_B\,\overline{F}(\ell_B,t_B),
\end{equation}
where $t_{\ell_A}$ is the time that particle $A$ reaches 
the site $\ell_A$ for the first time, and similar for 
$t_{\ell_B}$. The independence of the walkers is expressed by the fact 
that it is the product of two $\overline{F}-$functions which is being integrated.

Using \pref{j7} we have
\begin{equation}\label{j9}
P(t_{\ell_A}\le t_{\ell_B})=1-
 \int_0^\infty \,dt_A\,\overline{F}(\ell_A,t_A)\  
        \int_{0}^{t_A} \,dt_B\,\overline{F}(\ell_B,t_B)\,
\end{equation}
and interchanging the order of integration in this integral, we derive
the symmetry properties (which are evident from the random walk problem)
\begin{equation}\label{j10}
P(t_{\ell_A}\le t_{\ell_B})=1-P(t_{\ell_B}\le t_{\ell_A}),\quad 
{\rm hence} \quad P(t_{\ell_A}=t_{\ell_B})=\tfrac12, \quad 
{\rm if \ } \ell_A = \ell_B.
\end{equation}
If $\ell_A = \ell_B$ we can also use
integration by parts
\begin{equation}\label{j11}
           \begin{array}{l}
\dsp{\int_0^\infty \,dt_A\,\overline{F}(\ell,t_A)\  
        \int_{t_A}^\infty \,dt_B\,\overline{F}(\ell,t_B)=}\\
        \\
        \quad
          \dsp{-\int_0^\infty  
        \left[\int_{t_A}^\infty \overline{F}(\ell,\tau)\,d\tau\right]\,d\left[\int_{t_A}^\infty 
        \overline{F}(\ell,\tau)\,d\tau\right]=
        \tfrac12\left[\int_{0}^\infty \overline{F}(\ell,\tau)\,d\tau\right]^2=\tfrac12.}
\end{array}
\end{equation}

In this paper we derive asymptotic expansions of 
$P(t_{\ell_A}\le t_{\ell_B})$ given by
\begin{equation}\label{j12}
P= \ell_A\ell_B 
         \int_0^\infty \frac{dt_A}{t_A} e^{-t_A}I_{\ell_A}(t_A)\  
        \int_{t_A}^\infty \frac{dt_B}{t_B} e^{-t_B}I_{\ell_B}(t_B).
\end{equation}
In \pref{j12} the scale factor $\alpha$ has been absorbed in $t_A$ and $t_B$, because of 
$\frac{dt}{t}=\frac{d \alpha t}{\alpha t}$.

We will give one expansion that holds for large values of $\ell_A$ and one
for the case that both parameters $\ell_A$ and $\ell_B$ are large. We also
give an expansion that holds just when the sum $\ell_A+\ell_B$ is large.

\section{Transforming the integral}
\label{sec:trans}
We study the integral \pref{j12}.
We use well-known properties of the modified Bessel function to transform the double integral
in \pref{j12} into a single integral. 

The inner integral in \pref{j12} can be modified by evaluating
\begin{equation}\label{t1}
  S_{\ell}(t):= \ell\int_{t}^\infty  e^{-s}I_{\ell}(s)\frac{ds}{s}.
\end{equation}
where $\ell=1,2,\ldots\ $. We use the integral  representation (see
\cite[Eq. 9.6.19]{Abr})
\begin{equation}\label{t2}
I_n(s)=\frac1\pi\int_0^\pi e^{s \cos\theta} \cos n\theta \,d\theta.
\end{equation}
for integer values of $n$. Integrating by parts we obtain
\begin{equation}\label{t3}
\frac{\ell}{s} I_\ell(s)=\frac{1}\pi\int_0^\pi e^{s \cos\theta} \sin\theta\,\sin\ell\theta \,d\theta.
\end{equation}
It follows that
\begin{equation}\label{t4}
           \begin{array}{lll}
        \dsp{\ell\int_t^\infty e^{-ps} I_\ell(s)\,\frac{ds}{s}}&=&
        \dsp{\frac1\pi\int_0^\pi d\theta \sin\theta\,\sin\ell\theta \int_t^\infty ds\,e^{-s(p-\cos\theta)}}\\
          \\
          &=& \dsp{\frac{e^{-pt}}{\pi}\int_0^\pi\frac{\sin\theta\,\sin\ell\theta}{p-\cos\theta}\,e^{t\cos\theta}\,d\theta,}
          \end{array}
\end{equation}
which holds for $p\ge1$.
It follows that $S_{\ell}(t)$ of \pref{t1} can be written as
\begin{equation}\label{t5}
      S_{\ell}(t)= \frac{e^{-t}}{\pi}\int_0^\pi \frac{\sin\theta\,\sin\ell\theta}{1-\cos\theta}\,
       e^{t\cos\theta}\,d\theta=\frac{e^{-t}}{\pi}\int_0^\pi \cot\tfrac12\theta\, 
      \sin \ell\theta\,
      e^{t\cos\theta}\,d\theta.
\end{equation}
Using this relation and interchanging the order of integration in \pref{j12},
we obtain
\begin{equation}\label{t6}
 \begin{array}{l}
P= \dsp{\ell_A \int_0^\infty e^{-t}I_{\ell_A}(t) S_{\ell_B}(t)\frac{dt}{t}}\\ \\
        \quad = \dsp{\frac{\ell_A}{\pi}\int_0^\pi \cot\tfrac12\theta\, 
      \sin \ell_B\theta\,\left[\int_0^\infty
      e^{-2t+t\cos\theta}\,I_{\ell_A}(t) \frac{dt}{t}\right]d\theta}. 
\end{array}
\end{equation}
Invoking again \pref{t3} we obtain
\begin{equation}\label{t7}
P=\frac1{\pi^2}\int_0^\pi d\theta_2\,\frac{\sin\theta_2\,\sin\ell_B\theta_2}{1-\cos\theta_2}
\int_0^\pi d\theta_1\,\frac{\sin\theta_1\,\sin\ell_A\theta_1}{(1-\cos\theta_1)+(1-\cos\theta_2)}.
\end{equation}
The $\theta_1-$integral can be evaluated; see \cite[Eq. 3.613(3)]{Grad}. 
Another way is to use in \pref{t6} the Laplace integral 
\begin{equation}\label{t9}
\ell\int_0^\infty e^{-pt} I_\ell(t)\,\frac{dt}{t}=
\left(p+\sqrt{p^2-1}\right)^{-\ell}, \quad \ell>0, \quad p\ge 1,
\end{equation}
which follows from  \cite[29.3.53]{Abr} by taking $a=1, b=-1$. This 
gives
\begin{equation}\label{t10}
P=\frac{1}{\pi}\int_0^\pi \cot\tfrac12\theta\, 
      \sin \ell_B\theta\,\left(p+\sqrt{p^2-1}\right)^{-\ell_A}d\theta,
      \quad p=2-\cos\theta.
\end{equation}

\section{Asymptotic expansions}
\label{sec:as}
We give three asymptotic expansions: 
\begin{itemize}
\item
one for large $\ell_A$, with 
$\ell_B$ fixed, or small, 
\item
one for large $\ell_A$ and $\ell_B$, with $\ell_A\sim\ell_B$,
\item
one uniform expansion in which one or 
both parameters may be large.
\end{itemize}
\subsection{The case {\protect\boldmath $\ell_A\gg\ell_B$}}
\label{sec:as.a}
We start from \pref{t10} in the form
\begin{equation}\label{asa1}
P=\frac{1}{\pi}\int_0^\pi f(\theta) e^{-\ell_A\phi(\theta)}\,d\theta,
\end{equation}
where
\begin{equation}\label{asa2}
f(\theta)= \cot\tfrac12\theta\, 
      \sin \ell_B\theta, \quad 
       \phi(\theta)=\ln\left(p+\sqrt{p^2-1}\right),
      \quad p=2-\cos\theta.
\end{equation}
First we observe that
\begin{equation}\label{asa3}
\phi'(\theta)=\frac{\sin\theta}{\sqrt{p^2-1}}=\frac{\cos\frac12\theta}{\sqrt{1+\sin^2\frac12\theta}}.
\end{equation}
Hence, $\phi(\theta)$ is an increasing function on $[0,\pi]$ with
\begin{equation}\label{asa4}
\phi(0)=0, \quad \phi'(0)=1, \quad \phi'(\pi)=0.
\end{equation}
It follows that
\begin{equation}\label{asa5}
P\sim\frac{1}{\pi}\int_0^{\theta_0} f(\theta) e^{-\ell_A\phi(\theta)}\,d\theta,
\end{equation}
where $\theta_0$ is a fixed number in $(0,\pi)$, and the error in 
this approximation is exponentially small when $\ell_A$ is large.

Carrying out an integration by parts in the form
\begin{equation}\label{asa6}
P\sim\frac{-1}{\pi\ell_A}\int_0^{\theta_0} \frac{f(\theta)}{\phi'(\theta)}\,
d e^{-\ell_A\phi(\theta)}
\end{equation}
leads to
\begin{equation}\label{asa7}
P\sim\frac{-1}{\pi\ell_A} 
\frac{f(\theta)}{\phi'(\theta)}
e^{-\ell_A\phi(\theta)}\Biggr|_{\theta=0}^{\theta=\theta_0}+
\frac{1}{\pi\ell_A}\int_0^{\theta_0} f_1(\theta) e^{-\ell_A\phi(\theta)} \,d\theta,
\end{equation}
where
\begin{equation}\label{asa8}
f_1(\theta)= \frac{d}{d\theta}\frac{f(\theta)}{\phi'(\theta)}. 
\end{equation}
We can repeat this procedure, and compute the integrated terms. 
The terms at $\theta_0$ can be neglected because  
they give exponentially small contributions compared with 
the contributions from $\theta=0$. Note that we cannot 
take $\theta_0=\pi$, because $\phi'(\pi)=0$.

In this way we obtain the asymptotic expansion
\begin{equation}\label{asa9}
P\sim \frac{1}{\pi\ell_A}\left[a_0+\frac{a_1}{\ell_A}+
\frac{a_2}{\ell^2_A}+\ldots\right],
\end{equation}
where, for $k=0,1,2,\ldots$,
\begin{equation}\label{asa10}
a_k=\frac{f_k(0)}{\phi'(0)},\quad 
f_{k+1}(\theta)= \frac{d}{d\theta}\frac{f_{k}(\theta)}{\phi'(\theta)}, \quad
f_{0}(\theta)=f(\theta).
\end{equation}
The coefficients $a_k$ with odd indices are zero.
This follows from observing that $f(\theta)$ and $\phi'(\theta)$ are even functions;
see \pref{asa2} and \pref{asa3}.
Hence, $f_1(\theta)$ of \pref{asa8} is odd. By using the recursion in
\pref{asa10} it follows that $f_{2k}(\theta)$ is even, and that 
$f_{2k+1}(\theta)$ is odd. 
The first non-zero coefficients are
\begin{equation}\label{asa11}
a_0=2\ell_B, \quad 
a_2= \tfrac23 \ell_B(1-\ell^2_B),\quad
a_4= \tfrac1{30}\ell_B(23-80 \ell^2_B+12\ell^4_B).
\end{equation}

\begin{figure}
\begin{center}
\epsfxsize=8cm \epsfbox{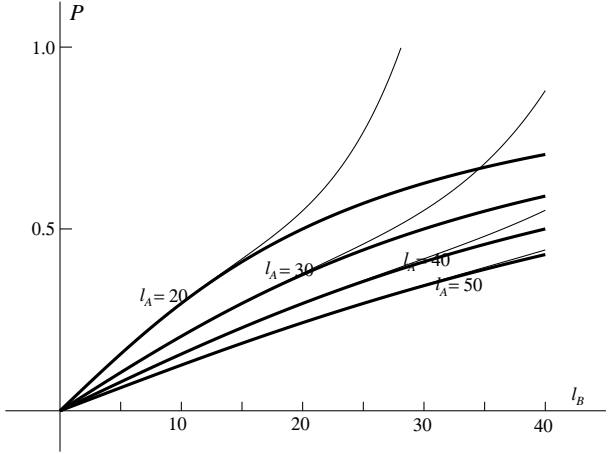}
\end{center}
\caption{\small Graphs of $P(t_{\ell_A}\le t_{\ell_B})$ based on the
asymptotic approximation \pref{asa9} (shown in red),  compared with graphs based on the expansion
\pref{asu13} (shown in blue). The red and blue graphs deviate from each other because
of the failure of the red approximations for large values of $\ell_B$. \label{walks2}
}
\end{figure}

In Figure \ref{walks2} we compare the approximations based on \pref{asa9} with values
obtained by using the expansion in \pref{asu13}, which holds when $\ell_A+\ell_B$
is large. We see that for smaller values of $\ell_B$ the graphs of the asymptotic
approximation \pref{asa9} are in agreement with the graphs obtained from the
expansion that holds when at least one of the parameters $\ell_A$ or $\ell_B$ is
large. The failure of the red approximations  is due to the failure of the asymptotic 
approximation \pref{asa9} that has been chosen for this case.

\subsection{The case {\protect\boldmath $\ell_A\sim\ell_B$}, both large}
\label{sec:as.aisb}
We replace in \pref{t10} $\ell_A$ by $\ell$ and $\ell_B$ by $\ell+\delta$. 
We know that $P=\frac12$ if $\delta=0$. 
We expand \pref{t10} for large values of $\ell$, 
keeping $\delta$ fixed. We have
\begin{equation}\label{sim1}
P=\tfrac12+P_1+P_2,
\end{equation}
where
\begin{equation}\label{sim2}
           \begin{array}{l}
        \dsp{P_1=\frac{1}{\pi}\int_0^\pi \cot\tfrac12\theta\, 
      (\cos\delta\theta-1)\,\sin \ell\theta 
      \,\left(p+\sqrt{p^2-1}\right)^{-\ell}d\theta}\\
      \\
      \quad\ 
      \dsp{=\Im\frac{1}{\pi}\int_0^\pi \cot\tfrac12\theta\, 
      (\cos\delta\theta-1)\,e^{ i\ell\theta} 
      \,\left(p+\sqrt{p^2-1}\right)^{-\ell}d\theta,}\\
      \\
      \dsp{      
P_2=\frac{1}{\pi}\int_0^\pi \cot\tfrac12\theta\, 
      \sin\delta\theta \,\cos \ell\theta \,
      \,\left(p+\sqrt{p^2-1}\right)^{-\ell}d\theta}\\
      \\
      \quad\ 
      \dsp{=\Re\frac{1}{\pi}\int_0^\pi \cot\tfrac12\theta\, 
      \sin\delta\theta \,e^{i\ell\theta} \,
      \,\left(p+\sqrt{p^2-1}\right)^{-\ell}d\theta.}
\end{array}
\end{equation}
Let
\begin{equation}\label{sim3}
g(\theta)=\cot\tfrac12\theta\,(\cos\delta\theta-1), \quad
\phi(\theta)=-i\theta+\ln\left(p+\sqrt{p^2-1}\right).
\end{equation}
Then we integrate by parts in the integral for $P_1$:
\begin{equation}\label{sim4}
\frac{1}{\pi}\int_0^\pi
g(\theta)e^{-\ell\phi(\theta)}\,d\theta=-
\frac{1}{\pi\ell}\int_0^\pi
\frac{g(\theta)}{\phi'(\theta)}\,de^{-\ell\phi(\theta)},
\end{equation}
and we obtain an expansion as in \pref{asa9},
\begin{equation}\label{sim5}
P_1\sim \frac{1}{\pi\ell}\left[b_0+\frac{b_1}{\ell}+
\frac{b_2}{\ell^2}+\ldots\right],
\end{equation}
where, for $k=0,1,2,\ldots$,
\begin{equation}\label{sim6}
b_k=\Im\frac{g_k(0)}{\phi'(0)},\quad 
g_{k+1}(\theta)= \frac{d}{d\theta}\frac{g_{k}(\theta)}{\phi'(\theta)}, \quad
g_{0}(\theta)=g(\theta).
\end{equation}
It turns out that the coefficients with even indices are zero. 
To verify this we can use a similar argument
as for the $a_k$ in \pref{asa10}.
The first non-zero coefficients are
\begin{equation}\label{sim7}
b_1=-\tfrac12\delta^2, \quad 
b_3= \tfrac14\delta^2,\quad
b_5= \tfrac{1}{96}\delta^2(4\delta^4-20\delta^2-77).
\end{equation}

In a similar way, let 
$h(\theta)=\cot\tfrac12\theta\,\sin\delta\theta$. Then 
\begin{equation}\label{sim8}
P_2\sim \frac{1}{\pi\ell}\left[c_0+\frac{c_1}{\ell}+
\frac{c_2}{\ell^2}+\ldots\right],
\end{equation}
where, for $k=0,1,2,\ldots$,
\begin{equation}\label{sim9}
c_k=\Re\frac{h_k(0)}{\phi'(0)},\quad 
h_{k+1}(\theta)= \frac{d}{d\theta}\frac{h_{k}(\theta)}{\phi'(\theta)}, \quad
h_{0}(\theta)=h(\theta).
\end{equation}
It turns out that the coefficients with even indices are zero.
The first non-zero coefficients are
\begin{equation}\label{sim10}
c_0= \delta, \quad 
c_2= \tfrac16\delta(\delta^2-1),\quad
c_4= -\tfrac{1}{240}\delta^2(12\delta^4+20\delta^2-77).
\end{equation}

\subsection{The case {\protect\boldmath $\ell_A+\ell_B$} large}
\label{sec:as.u}
Because of (see \pref{asa4})
\begin{equation}\label{asu1}
e^{-\ell_A\phi(\theta)}\sim e^{-\ell_A\theta},\quad \theta\to 0,
\end{equation}
we have for large values of $\ell_A$
\begin{equation}\label{asu2}
P\sim\frac{2}{\pi}\int_0^\infty
      \frac{\sin \ell_B\theta}{\theta}\,e^{-\ell_A\theta}\, d\theta=
      \frac2\pi\,\arctan\frac{\ell_B}{\ell_A},
\end{equation}
where we used  \cite[Eq. 29.3.110]{Abr}.

Observe that this estimate perfectly reflects the properties of $P$ 
mentioned in \pref{j10}; also, it is less than unity, as the probability 
$P$ itself is.  Moreover, in this estimate large values of $\ell_B$ 
do not disturb the approximation.

The result \pref{asu2} is obtained by combining the dominant behaviour of 
$e^{-\ell_A\phi(\theta)}$ near the origin with the complete form 
$\sin \ell_B\theta$, without expanding this function. 

We modify the integration by parts procedure of \S \ref{sec:as.a},
by including the (possible large) parameter $\ell_B$ in 
the "phase function"  $\phi(\theta)$. We can do this by writing 
$\sin \ell_B\theta=\Im e^{i\ell_B\theta}$.
A complication is the pole of the function 
$\cot\tfrac12\theta$, which singularity is removable in combination with
the function $\sin \ell_B\theta$.

To perform the integration by parts procedure we proceed in the following way.
In \pref{t10} we can consider $\ell_B$ as a continuous parameter, and we 
can
differentiate with respect to $\ell_B$. We also observe that $P$ vanishes 
with $\ell_B$. We have
\begin{equation}\label{asu3}
\frac{\partial P}{\partial \ell_B}=\frac{1}{\pi}\Re\left[\int_0^\pi 
\theta \cot\tfrac12\theta\, 
      e^{i\ell_B\theta}\,\left(p+\sqrt{p^2-1}\right)^{-\ell_A}d\theta\right].
\end{equation}
We write this in the form
\begin{equation}\label{asu4}
\frac{\partial P}{\partial \ell_B}=\frac{2}{\pi}\Re\, Q,
\end{equation}
where
\begin{equation}\label{asu5}
Q= \int_0^\pi 
f(\theta)e^{\psi(\theta)}\,d\theta,
\end{equation}
with
\begin{equation}\label{asu6}
f(\theta)=\tfrac12\theta\cot\tfrac12\theta,\quad \psi(\theta)=i\ell_B\theta-\ell_A
\ln(p+\sqrt{p^2-1}).
\end{equation}
We integrate by parts, starting with
\begin{equation}\label{asu7}
Q= \int_0^\pi 
\frac{f(\theta)}{\psi'(\theta)}\,de^{\psi(\theta)}=
\frac{f(\theta)}{\psi'(\theta)}\,e^{\psi(\theta)}\Biggr|_{\theta=0}^{\theta=\pi}
+\int_0^\pi 
f_1(\theta)e^{\psi(\theta)}\,d\theta,
\end{equation}
where
\begin{equation}\label{asu8}
f_1(\theta)=-\frac{d}{d\theta}\frac{f(\theta)}{\psi'(\theta)}.
\end{equation}
We repeat this procedure, and compute the integrated terms. 
Again, the terms at $\theta=\pi$ can be neglected.

We obtain
\begin{equation}\label{asu9}
Q\sim d_0+d_1+d_2+\ldots,
\end{equation}
where
\begin{equation}\label{asu10}
d_k=-\frac{f_k(0)}{\psi'(0)}, \quad 
f_k(\theta)=-\frac{d}{d\theta}\frac{f_{k-1}(\theta)}{\psi'(\theta)}, \quad
k=0,1,2,\ldots,
\end{equation}
and $f_0(\theta)=f(\theta)$. Again, all coefficients with odd index vanish. 
This follows from
\begin{equation}\label{asu11}
\psi'(\theta)=i\ell_B-\ell_A\frac{\sin\theta}{\sqrt{p^2-1}}=
i\ell_B-\ell_A\frac{\cos\frac12\theta}{\sqrt{1+\sin^2\frac12\theta}},
\end{equation}
which is an even function and $f(\theta)$ is also even. Hence,
$f_1(\theta)$ in \pref{asu8} is odd; and so on.

We have
\begin{equation}\label{asu12}
d_0=\frac1{\ell_A-i\ell_B},\quad 
d_2=\frac{2\ell_A+i\ell_B}{6(\ell_A-i\ell_B)^4}, \quad
d_4=\frac{23\ell^2_A+129 i\ell_A \ell_B+2\ell^2_B}{60(\ell_A-i \ell_B)^7}.
\end{equation}

Considering \pref{asu4},
taking the real parts of 
the coefficients and integrating the real parts 
over the interval $[0,\ell_B]$, we find    
\begin{equation}\label{asu13}
P\sim\frac2\pi\left(e_0+e_2+e_4+e_6\ldots\right),
\end{equation}
where
\begin{equation}\label{asu14}
e_{2k}=\int_0^{\ell_B} d_{2k}(\ell'_B)\, d\ell'_B.
\end{equation}
The first few are
{\small
\begin{equation}\label{asu15}
\begin{array}{l}
e_0=\dsp{\arctan\frac{\ell_B}{\ell_A}}, \\ \\ 
e_2=\dsp{\frac{\ell_A\ell_B(\ell_A^2-\ell_B^2)}{3(\ell_A^2+\ell_B^2)^3}},\\ \\
e_4=\dsp{\frac{\ell_A\ell_B(\ell_A^2-\ell_B^2)
(23\ell_A^4-354\ell_A^2\ell_B^2+23\ell_B^4)}
{60(\ell_A^2+\ell_B^2)^6}},\\ \\
e_6=\dsp{\frac{\ell_A\ell_B(\ell_A^2-\ell_B^2)
(249\ell_A^8-10796\ell_A^6\ell_B^2+40630\ell_B^4\ell_A^4-10796\ell_B^6\ell_A^2+249\ell_B^8)}
{126(\ell_A^2+\ell_B^2)^9}}.
\end{array}
\end{equation}
}%
We see that the shown coefficients $e_2, e_4, e_6$ vanish when 
$\ell_A=\ell_B$, and that in fact $e_{2k}(\ell_A,\ell_B)=\frac12\pi\delta_{k,0}-e_{2k}(\ell_B,\ell_A)$,
$k=0,1,2,\ldots$.
These properties are in agreement with the relations for $P$ in \pref{j10}.
Because there is no symmetry in \pref{t10} with respect to $\ell_A$ 
and $\ell_B$, they do not follow from the construction of the coefficients 
$d_{2k}$ and $e_{2k}$. 

When we scale the parameters by putting $\ell_B=\lambda \ell_A$, we see 
that the shown coefficients obey the relation
\begin{equation}\label{asu16}
d_{2k}=\bo\left(\ell_A^{-2k}\right),
\end{equation}
uniformly with respect to $\lambda\ge0$. When we write \pref{asu9} with a 
remainder, that is,
\begin{equation}\label{asu17}
Q= d_0+d_2+\ldots+d_{2k-2}+\int_0^\pi f_{2k}(\theta)e^{\psi(\theta)}\,d\theta,
\end{equation}
a straightforward analysis shows that similarly
\begin{equation}\label{asu18}
f_{2k}(\theta)=\bo\left(\ell_A^{-2k}\right),
\end{equation}
uniformly with respect to $\lambda\ge0$ and $\theta\in[0,\pi]$.
This shows the nature of the uniform asymptotic expansion of Q, and, after integrating, 
the nature of the expansion for the probability $P$.

By expanding the coefficients $e_k$ in \pref{asu15} for large $\ell_A$ 
with $\ell_A\gg\ell_B$, we obtain the coefficients of the non-uniform
expansion of \S\ref{sec:as.a}.

\section{Three particles and more}\label{sec:three}
For three random walkers $A, B, C$ the probability integral reads
\begin{equation}\label{d1}
           \begin{array}{l}
\dsp{P(t_{\ell_A}\le t_{\ell_B}\le t_{\ell_C})=}\\
\\
\quad\quad\quad \dsp{  \int_0^\infty \,dt_A\,\overline{F}(\ell_A,t_A)\
        \int_{t_A}^\infty \,dt_B\,\overline{F}(\ell_B,t_B)\
\int_{t_B}^\infty \,dt_C\,\overline{F}(\ell_C,t_C),}
\end{array}
\end{equation}
with the density as in \pref{j8}. That is,
\begin{equation}\label{d2}
           \begin{array}{l}
\dsp{P(t_{\ell_A}\le t_{\ell_B}\le t_{\ell_C})=}\\
\\
\ \dsp{ \ell_A\ell_B\ell_C \int_0^\infty \,\frac{dt_A}{t_A}\,e^{-t_A}I_{\ell_A}(t_A)\
        \int_{t_A}^\infty \,\frac{dt_B}{t_B}\,e^{-t_B}I_{\ell_B}(t_B)\
\int_{t_B}^\infty \,\frac{dt_C}{t_C}\,e^{-t_C}I_{\ell_C}(t_C).}
\end{array}
\end{equation}
It gives the probability that particle $A$ reaches site
$\ell_A$, before particle $B$ reaches $\ell_B$, while
$B$ reaches site
$\ell_B$, before particle $C$ reaches $\ell_C$.

First we observe that the probability for
three particles arriving at the same  site $\ell$, that is,  
$\ell_A=\ell_B=\ell_C=\ell$ 
equals $\frac1{3!}$. This easily follows from (cf. \pref{j11})
\begin{equation}\label{d3}
 \int_{t_A}^\infty \,dt_B\,\overline{F}(\ell,t_B)\ 
\int_{t_B}^\infty \,dt_C\,\overline{F}(\ell,t_C)=
\tfrac12\left[\int_{t_A}^\infty \overline{F}(\ell,\tau)\,d\tau\right]^2.
\end{equation}
Substituting this in \pref{d2},  performing another integration by parts,
and using \pref{j7}, gives the value $\frac1{3!}$. 
Using the same method we infer that for $n$ particles
the probability for all $n$ particles arriving at the same site $\ell$
equals $\frac1{n!}$.

Repeating the steps used for obtaining \pref{t7}, and replacing all Bessel functions 
by using \pref{t3}, we easily find for \pref{d3}
\begin{equation}\label{d4}
           \begin{array}{l}
\dsp{P(t_{\ell_A}\le t_{\ell_B}\le t_{\ell_C})=
\frac1{\pi^3}\int_0^\pi\,d\theta\,\frac{\sin\theta\,\sin\ell_C\theta}{1-\cos\theta}\ \times}\\
\\
\quad\quad \dsp{ 
\int_0^\pi\, d\sigma\,\frac{\sin\sigma\,\sin\ell_B\sigma}{2-\cos\theta-\cos\sigma}\int_0^\pi
\,d\tau\,\frac{\sin\tau\,\sin\ell_A\tau}{3-\cos\theta-\cos\sigma-\cos\tau}.}
\end{array}
\end{equation}
Evaluating the $\tau-$integral gives
\begin{equation}\label{d5}
P=\frac{1}{\pi^2}\int_0^\pi\,\int_0^\pi \, 
\frac{\sin\theta\,\sin\ell_C\theta}{1-\cos\theta}
\frac{\sin\sigma\,\sin\ell_B\sigma}{2-\cos\theta-\cos\sigma}
 \left(q+\sqrt{q^2-1}\right)^{-\ell_A}\,d\theta\,d\sigma,   
\end{equation}
where $q=3-\cos\theta-\cos\sigma.$

From the above analysis it is clear how a similar integral representation can 
be obtained for $n$ random walkers $A_1, A_2, \ldots, A_n$.
The probability can be written in the form of the $n-$fold integral
\begin{equation}\label{d6}
P(t_{\ell_{A_1}}\le t_{\ell_{A_2}}\le \ldots\le t_{\ell_{A_n}})=
\frac{1}{\pi^{n}}\int_0^\pi\,d\theta_1\cdots \int_0^\pi\,d\theta_n\,
         \prod_{j=1}^{n}\frac{\sin\theta_j\,\sin\ell_{A_j}\theta_j}{\wt{p}_j},
\end{equation}
where
\begin{equation}\label{d7}
\wt{p}_j=\sum_{k=j}^n(1-\cos\theta_k)=2\sum_{k=j}^n\sin^2\tfrac12\theta_k,\quad j=1,2,\ldots, n.
\end{equation}
Integrating the $\theta_1$ integral  gives
\begin{equation}\label{d8}
P=
\frac{1}{\pi^{n-1}}\int_0^\pi\,d\theta_2\cdots \int_0^\pi\,d\theta_n\,\left(p+\sqrt{p^2-1}\right)^{-\ell_{A_1}}\,
         \prod_{j=2}^{n}\frac{\sin\theta_j\,\sin\ell_{A_j}\theta_j}{\wt{p}_j},
\end{equation}
where
\begin{equation}\label{d9}
p=1+\wt{p}_2=n-\sum_{j=2}^n\cos\theta_j.
\end{equation}

\subsection{Asymptotic approximations for three particles}
For large values of $\ell_A$ the main contributions to the integral in \pref{d5} 
come from the origin $\sigma=0$, $\theta=0$. To see this we observe that
\begin{equation}\label{d10}
           \begin{array}{l}
\dsp{q+\sqrt{q^2-1}=}\\
\\
\quad
                \dsp{3-\cos\theta-\cos\sigma+\sqrt{(2-\cos\theta-\cos\sigma)(4-\cos\theta-\cos\sigma)}=}\\
                \\
\quad
                \dsp{1+\sqrt{\theta^2+\sigma^2}+\bo\left(\theta^2,\theta\sigma,\sigma^2\right)},
      \end{array}
\end{equation}
and that
\begin{equation}\label{d11}
\left(q+\sqrt{q^2-1}\right)^{-\ell_A}=e^{-\ell_A \ln(q+\sqrt{q^2-1})}\sim 
e^{-\ell_A\sqrt{\theta^2+\sigma^2}},
\end{equation}
as $\theta, \sigma \to 0$. We also have
\begin{equation}\label{d12}
\frac{\theta\sin\theta}{1-\cos\theta}\sim2,\quad
\frac{\sigma\sin\sigma}{2-\cos\theta-\cos\sigma}\sim\frac{2\sigma^2}{\theta^2+\sigma^2}
\end{equation}
as $\theta, \sigma \to 0$.

This motivates us to consider as a first approximation
\begin{equation}\label{d13}
P\sim \frac {4}{\pi^2} \int_0^\pi\,d\theta\,\int_0^{\pi} \,d\sigma\,
\frac{\sin(\ell_C\theta)}{\theta}\, 
\frac{\sin(\ell_B\sigma)}{\sigma}\,\frac{\sigma^2}{\theta^2+\sigma^2} 
e^{-\ell_A\sqrt{\theta^2+\sigma^2}},
\end{equation}
where we have used \pref{d5}, \pref{d11} and \pref{d12}.

Next we use polar coordinates for $\theta$ and $\sigma$ by writing
\begin{equation}\label{d14}
\theta = r \cos\phi,\quad  \sigma = r \sin\phi, 
\quad 0\le r\le \pi,\quad 0\le\phi\le\tfrac12\pi.
\end{equation}
We extend the finite square in the $(\theta,\sigma)-$plane to the quarter plane and obtain
\begin{equation}\label{d15}
P\sim \frac {4}{\pi^2} \int_0^{\frac12\pi}\,d\phi\,\int_0^{\infty} r\,dr\,
\frac{\sin(\ell_C r\cos\phi)}{r\cos\phi}\, 
\frac{\sin(\ell_B r\sin\phi)}{r\sin\phi}\, 
\sin^2\phi\, e^{-\ell_Ar} \,.
\end{equation}
The $r$ integral can be found in \cite[Eq. (3.947)]{Grad}, that is,
\begin{equation}\label{d16}
\int_0^{\infty} 
e^{-ar} \sin (br) \sin (cr)  \,\frac{dr}{r}=
\tfrac14\ln\frac{a^2+(b+c)^2}{a^2+(b-c)^2},
\end{equation}
and can be proved by differentiation with respect to $a$. We obtain
\begin{equation}\label{d17}
P\sim \frac {1}{\pi^2} \int_0^{\frac12\pi}\tan\phi 
\ln\frac{\ell_A^2+(\ell_C\cos\phi+\ell_B\sin\phi)^2}{\ell_A^2+(\ell_C\cos\phi-\ell_B\sin\phi)^2}\,d\phi,
\end{equation}
which can be written as
\begin{equation}\label{d18}
P\sim \frac {1}{\pi^2} \int_0^{\frac12\pi}\tan\phi 
\ln\frac{1+u\cos2\phi+v\sin2\phi}{1+u\cos2\phi-v\sin2\phi}\,d\phi,
\end{equation}
where
\begin{equation}\label{d19}
u= \frac{\ell_C^2-\ell_B^2}{2\ell_A^2+\ell_B^2+\ell_C^2},\quad v=\frac{2\ell_C \ell_B}{2\ell_A^2+\ell_B^2+\ell_C^2}.
\end{equation}
When $v$ is small we can expand
\begin{equation}\label{d20}
\begin{array}{rcl}
\dsp{\ln\frac{1+u\cos2\phi+v\sin2\phi}{1+u\cos2\phi-v\sin2\phi}}
&=&
\dsp{\ln\frac{1+\frac{v\sin2\phi}{1+u\cos2\phi}}{1-\frac{v\sin2\phi}{1+u\cos2\phi}}}\\
& &\\
&=&2
\dsp{\sum_{n=0}^\infty
 \frac1{2n+1} \frac{v^{2n+1}\sin^{2n+1}2\phi}{(1+u\cos2\phi)^{2n+1}}},
\end{array}
\end{equation}
which gives
\begin{equation}\label{d21}
P\sim\frac2{\pi^2}\sum_{n=0}^\infty \frac{v^{2n+1}}{2n+1} \int_0^{\frac12\pi}
\tan \phi \frac{\sin^{2n+1}2\phi}{(1+u\cos2\phi)^{2n+1}}\,d\phi.
\end{equation}
This expansion is useful when $\ell_A$ is large compared with $\ell_B$ and $\ell_C$.

The  integral in \pref{d21} can be written in terms of a 
Gauss hypergeometric function,
and the sum can be written as an Appell function. This does not give further insight, however.
We prefer to give a few further estimates.

For examining the convergence of the series in \pref{d21}, observe that
\begin{equation}\label{d22}
\tan \phi \frac{\sin^{2n+1}2\phi}{(1+u\cos2\phi)^{2n+1}}\le 
\frac{2}{(1-u)^{2n+1}},
\end{equation}
with
\begin{equation}\label{d23}
1-u=\frac{2\ell_A^2+2\ell_B^2}{2\ell_A^2+\ell_B^2+\ell_C^2},
\end{equation}
which is bounded away from 0, unless $\ell_C$ is 
much larger than $\ell_A$ and $\ell_B$.

It follows that expansion \pref{d21}
can be viewed as an asymptotic expansion for small values of $v$
for the right-hand side in \pref{d18}.

Of further interest is that 
when $u=0$, that is, $\ell_B=\ell_C$ we can evaluate the right-hand side of
\pref{d21} in terms of elementary functions. In fact we obtain by using
\begin{equation}\label{d24}
\int_0^{\frac12\pi} \sin^{2n+2}\phi \cos^{2n}\phi\,d\phi=
\frac{\Gamma(n+\frac32)\Gamma(n+\frac12)}{\Gamma(2n+2)},
\end{equation}
a Gauss hypergeometric function, that can be written as an elementary function:
\begin{equation}\label{d25}
P\sim\frac1\pi v\, {}_2F_1\left(\tfrac12,\tfrac12;\tfrac32;v^2\right)
=\frac1\pi\arcsin v, \quad v=\frac{\ell_B^2}{\ell_A^2+\ell_B^2}.
\end{equation}
When $\ell_A=\ell_B=\ell_C$ this becomes
\begin{equation}\label{d26}
P\sim\frac1\pi \arcsin\frac12=\frac1{6}, 
\end{equation}
which is the exact value.

\section{Discussion and concluding remarks}\label{disc}
We have discussed in this paper a method  of considering different simultaneous
independent $1D-$random walks. This work has been motivated  by an attempt to describe the 
agglomeration of a number of random walkers on a linear chain which will be fixed
when they come to occupy nearest neighbour positions on the chain.
In treating this problem it turns out to be possible to effect a transformation of coordinates
which makes the evolution equation  become separable, such that we obtain 
a product of "one-particle" equations 
which can then be mathematically treated as independent 
random walkers as described in this paper. However, it turns out that this 
separation is possible only when the jump probabilities in both directions are equal. 
This is the reason why we have limited ourselves to equal jump probabilities 
in this work.

For two particles we have given a complete asymptotic description 
for the case when $\ell_A$ and/or
$\ell_B$ are large. For three particles we have also given asymptotic results,
but a full description becomes a very complicated matter.

Very recently a paper \cite{Dankel} has appeared which treats  a related problem 
(with discrete time steps)
by a different method, involving stochastic matrices.

\bigskip
\noindent
{\bf Acknowledgments\quad}\\
The work of the first author has been made possible by the kind 
hospitality of the AMOLF-Institute
of the Foundation FOM in Utrecht.\\
\noindent
N.M. Temme acknowledges financial support from Ministerio Ciencia y Tecnolog\'\i a
from project SAB2003-0113.\\
\noindent
The authors thank the referee for helpful suggestions to obtain the 
\lq angular\rq\ integral representations in a simpler way, in particular the
$n-$dimensional forms, and for several other improvements 
in the paper.

\end{document}